\documentclass[a4paper]{article}
\usepackage{latexsym}
\usepackage{amssymb}
\usepackage{amsmath}
\newtheorem{conj}{Conjecture}[section]
\newtheorem{lemma}[conj]{Lemma}
\newtheorem{coro}[conj]{Corollary} 

\newtheorem{prop}[conj]{Proposition}

\newtheorem{thm}{Theorem}
\newtheorem{rem}[conj]{Remark}

\newcommand{\qed}{\raisebox{-.8ex}{$\Box$}}
\newenvironment{bew}
{\noindent{\bf Proof.}}
{\hfill \qed\\}

\newcommand{\la}{\langle}
\newcommand{\ra}{\rangle}

\newcommand{\Aut}{{\rm Aut}}

\newcommand{\Syl}{{\rm Syl}}
\newcommand{\PGL}{{\rm PGL}}
\newcommand{\PSL}{{\rm PSL}}

\newcommand{\SL}{{\rm SL}}

\newcommand{\oF}{\overline{F}}
\newcommand{\oG}{\overline{G}}

\newcommand{\core} {{\rm core}}
\newcommand{\Stab} {{\rm Stab}}

\newcommand{\Sym}{~{\rm Sym}}

\newcommand{\Alt}{~{\rm Alt}}
\newcommand{\RMult}{~{\rm RMult}}

\newcommand{\ZZ}{\mathbb{Z}}

\newcommand{\QR}{{\rm I}\kern-5.0pt {\rm Q} \kern2pt}

\title{On Bruck Loops of 2-power Exponent\thanks{This research is part of the project ``Transversals in Groups with an application to loops'' GZ: BA 2200/2-2 funded by
the DFG}}
\author{B.Baumeister, G.Stroth, A.Stein}

\begin{document}

\maketitle
\begin{abstract}
The goal of this paper is two-fold. First we provide the information needed to study Bol, $A_r$ or Bruck loops
by applying group theoretic methods. This information is used in this paper as well as in [BS3] and in [S]. 

Moreover, we determine the groups  associated to Bruck loops of 2-power exponent under the assumption
that every nonabelian simple group $S$ is either passive or isomorphic to $\PSL_2(q)$, $q-1 \ge 4$ a $2$-power.
In a separate paper it is proven that indeed every nonabelian simple  group $S$ is either passive or isomorphic to $\PSL_2(q)$, $q-1 \ge 4$ a $2$-power [S].
The results obtained here are used in [BS3], where we determine the structure of the groups  associated to the
finite Bruck loops.
\end{abstract}

\section{Introduction}
For a long time the following groups $G$ have been studied:
\medskip\\
{\bf Hypothesis (A)}
{\sl Assume that $G$ has a subgroup $H$ such that there is a transversal $K$ to $H$ in $G$ which
is the union of $1 \in G$ and $G$-conjugacy classes of involutions.}
\medskip\\
It has been conjectured that $G$ is a $2$-group if $G$ is a finite group which is generated by $K$.

Nagy [Nag] as well as Baumeister and Stein [BS1] found a counterexample to that conjecture.
This paper is part of a series of papers where we  determine the structure of the finite groups appearing in those
triples $(G,H,K)$  which satisfy
Hypothesis (A). In the present paper we reduce the question on the structure of the groups  to a question on finite simple groups, which is solved in [BS2] and [S].

Notice that, as $K$ is closed under conjugation, $K$ is a transversal to all the conjugates
of $H$ in $G$, and moreover, $1 \in K$. Baer observed that we can construct out of such a triple
$(G,H,K)$ a loop [Baer] (see Section~2.1). A {\em loop} is a set $X$ together with a binary operation $\circ$ on $X$,
such that there exists a unique $1_{\circ} \in X$ with $1_{\circ} \circ x = x \circ 1_{\circ} = x$ for all $x \in X$ and
such that the left and right translations
\[ \lambda_x: X \to X, \quad y \mapsto x \circ y , \quad \rho_x: X \to X,\quad  y \mapsto y \circ x  \]
are bijections. A loop can be thought of as a non-associative group. 

Conversely given a loop $X$, we can recover the triple $(G,H,K)$ [Baer] (see Section~2.1). 
A triple $(G,H,K)$ is called {\em loop folder}, if
\begin{itemize}
 \item $K$ is a transversal to all the conjugates of $H$ and if
\item $1 \in K$.
\end{itemize}

Clearly, every triple $(G,H,K)$ satisfying (A) is a loop folder. Moreover, if (A) holds, then
$K$ is a {\em twisted subgroup} of $G$, that is $1 \in K$ and $ x^{-1}, xyx \in K$ for all
$x,y \in K$. This translates into the language of loops to the {\em right Bol Identity}:
\[ ((x \circ y) \circ z) \circ y = x \circ ((y \circ z) \circ y)  \quad \mbox{ for all} ~ x,y,z \in X, \]
where $(X, \circ)$ is the loop constructed from $(G,H,K)$.

A loop is called {\em (right) Bol loop}, if it satisfies this identity.
 In a Bol loop, the subloop generated by a single element is a cyclic group.
Therefore powers, inverses and orders of elements are well defined, as is the exponent of a finite Bol loop. 

The loop associated to a triple $(G,H,K)$ fulfilling (A) is a Bol loop of exponent $2$,
as $k^2 = 1 $ for every $k \in K$. The loop then also
satisfies the {\em  automorphic inverse property, AIP}, that is:
\[ (x \circ y)^{-1} = x^{-1} \circ y^{-1} \quad \mbox{ for all } x,y  \in X. \]

Bol loops with that property are {\em  Bruck loops}. Our project is not only to determine the structure of the groups in
the triples satisfying Hypothesis (A), but as well the structure of the groups appearing in the larger class of triples associated to the finite Bruck loops. In [BS3] we use the results proved in this paper and [S] to find out the structures of the 
possible $G, H$ and $K$.

In 2005 Aschbacher, Kinyon and Phillips gave insight into the structure of general finite  Bruck loops, as they showed in \cite{AKP}:
\begin{itemize}
\item Elements of 2-power order and elements of odd order commute.
\item Bruck loops are a central product of a subloop of odd order and a subloop generated by elements of 2-power order.
\item Simple Bruck loops are of 2-power exponent.
\item The structure of a minimal simple Bruck loop ($M$-loop) is very restricted (see Theorem~\ref{AMT}).
\end{itemize}
This focuses attention on {\em  Bruck loops of 2-power exponent}, i.e. Bruck loops where every element is of 2-power order.
We call a loop folder associated to a Bruck loop a {\em BX2P-folder}, if 
\begin{itemize}
 \item $K$ is a twisted subgroup
\item every element in $K$ has a $2$-power order
\item $H$ acts on $K$.
\end{itemize}

To formulate the  statement of the main theorem we need a further definition.
A finite nonabelian simple group $S$ is called {\em passive}, if whenever $(G,H,K)$ is a BX2P-folder with 
\[ F^\ast(G/O_2(G)) \cong S, \]
then $G = O_2(G) H$.

Notice that in this case the loop to $(G,H,K)$ is of 2-power size and therefore soluble  by \ref{solubleloop}.

\begin{thm}
\label{EnvelopeGroups}
Let $(G,H,K)$ be a loop folder associated to a finite Bruck loop such that:
\begin{itemize}
 \item[(a)] $K$ is a twisted subgroup
\item[(b)] every element in $K$ has a $2$-power order
\item[(c)] $H$ acts on $K$
\item[(d)] $G = \langle K \rangle$
\end{itemize}
 and assume, that every non-abelian simple section
of $G$ is either passive or isomorphic to $\PSL_2(q)$ for $q=9$ or a Fermat prime $q \ge 5$. Then the following holds: 
\begin{enumerate}
\item[(1)] $G/O_2(G) \cong D_1 \times D_2 \times \cdots \times D_e$ for some nonnegative integer $e$.
\item[(2)]  $D_i \cong \PGL_2(q_i)$ with $q_i \ge 5$ a Fermat prime or $q_i=9$, for $1 \le i \le e$.
\item[(3)] $D_i \cap \overline{H} \cong q_i:(q_i-1)$ is a Borel subgroup in $D_i$ with $\overline{H} := HO_2(G)/O_2(G)$.
\item[(4)] $F^\ast(G) = O_2(G)$
\end{enumerate}
\end{thm}

Notice that the assumption that the simple sections of $G$ are either passive or one of Aschbachers 
candidates is similiar to a ${\cal K}$-group assumption in the classification of finite simple groups, see Section 4.
Another way to think of the main theorem is as a structure reduction: 

Given any finite group $G$, are there $H$ and $K$, such that $(G,H,K)$ is a nice folder (i.e. (a) - (d) of Theorem~\ref{EnvelopeGroups} holds) to a Bruck loop?
The main theorem reduces this problem to the case of those $G$ such that $F^\ast(G/O_2(G))$ is a finite simple group.

We call a loop folder $(G,H,K)$ {\em nice} with respect to some loop property (Bol, $A_r$, Bruck) if this 
property translates into a group theoretic property of the triple $(G,H,K)$.

The only known example of a non-passive group is $\PSL_2(5)$. The work of Aschbacher, Kinyon and Phillips suggests that
 also $\PSL_2(q)$ is a non-passive group for other values of $q$ with $q-1$ a $2$-power.
Unfortunately it is an open question, whether these groups are passive or not. An answer demands either an example
or a proof of the nonexistence of examples. This relates to hard questions about 2-groups.
However in a forthcoming paper we show, that the non-passive finite simple groups are among the $\PSL_2(q)$ with 
$q-1 \ge 4$ a 2-power [S].

The structure of the finite groups $G$ which satisfy Hypothesis (A) and which are generated by $K$
is completely determined in [BS3]. Application of Theorem~\ref{EnvelopeGroups} yields the following:

\begin{coro}
Let $G$ be a finite group and $H \le G$, such that there is a transversal $K$ to $H$ in $G$ which
is the union of $1 \in G$ and $G$-conjugacy classes of involutions. If $G= \langle K \rangle$
and if  every non-abelian simple section of $G$ is either passive or isomorphic to $\PSL_2(q)$ for $q=9$ or a Fermat prime $q \ge 5$, then $(G,H,K)$
is a loop envelope to a Bruck loop of  exponent $2$ with $H$ acting on $K$. Therefore, Theorem~\ref{EnvelopeGroups}
 describes $G$, $H$ and $K$.
\end{coro}

The organisation of the paper is as follows: In the next section we introduce the relevant notation on loops and
 assemble all the important facts on the relation between loops and groups. The idea of that section is to
provide a base for our series of paper - the results given there will be needed in this paper, in [S] and in [BS3] -
as well as for future papers on loops.
If we have a proof of some result which is more instructive then the known one, then we include that proof. 
Else we quote the literature.
 The third section contains general results on Bruck loops of $2$-power exponent,
which provide a set of tools for the classification of non-passive groups. Finally, the proof
of the main theorem is contained in Section~4.

\section{The relation between loops and groups}

We follow the notation of Aschbacher \cite{A} and \cite{AKP}. In particular we use the right Bol identity and
talk about right Bol loops. As there is an opposite relation between left and right Bol loops,
the decision between left and right Bol loops is only a notational convention, but also in the tradition of Bol, Bruck, Glauberman
and Aschbacher.

\subsection{The Baer correspondence}

Baer observed that statements about loops can be translated into the language of group theory \cite{B}.

Given a loop $(X,\circ)$, we define for $x \in X$ a map $\rho: X \to Sym(X), x \mapsto \rho_x$.
 
We record some standard  loop theoric notation.
\begin{itemize}
\item[] $G := \RMult(X) := \langle \rho_x: x \in X\rangle \le Sym(X)$, 
\item[] $H:=\Stab_G(1_{\circ})$,
\item[] $K:=\{ \rho_x: x \in X \} \subseteq G$ and 
\item[] $\kappa: K \to X : \rho_x \mapsto x $.  
\end{itemize}

Then $(G,H,K)$ satisfies the following properties:
\begin{enumerate}
\item[(1)] $K$ is a transversal to all conjugates of $H$.
\item[(2)] $H$ is core free: $1 = \cap_{g \in G} H^g$.
\item[(3)] $G = \langle K \rangle $.
\item[(4)] $1 \in K$. 
\end{enumerate}

\noindent
{\bf Definition}\label{folder_defi}
A triple $(G,H,K)$ with $G$ a group, $H \le G$ and $K \subseteq G$ is called 
\begin{itemize}
\item[] a {\em loop folder}, if it satisfies (1) and (4),
\item[] a {\em faithful loop folder}, if it satisfies (1) and (2),
\item[] a {\em loop envelope}, if it satisfies (1), (3) and (4),
\end{itemize}

\noindent
{\bf Remarks} 
(a) (1) is equivalent to the property\\
(1'): \hspace{2cm} $|K|=|G:H|~\mbox{and}~ H^g \cap K K^{-1} =1~\mbox{for all}~g \in G.$

(b) (1) and (2) imply (4). 

(c) Conditions (2) and (3) seem to be natural, but may not be satisfied in
loop folders to subloops, so called subfolders (see below for a definition).

(d) To a loop $(X,\circ)$, there is up to isomorphism (of loop folders) a unique loop
folder to $X$ satisfying (2) and (3): The loop folder, which we constructed above in $G=\RMult(X)$. We call it the {\em Baer envelope} of the loop.

(e) Let $(G,H,K)$ be a loop folder and let $\kappa$ be a bijection between $K$ and some set $X$.
Then the following operation $\circ$ on $X \times X$ defines a loop on $X$:
Set for all $k_1,k_2 \in K$
 $$\kappa(k_1) \circ \kappa(k_2) = \kappa(k_{12})~\mbox{where}~\{ k_{12} \} = K \cap H k_1 k_2.$$
(Notice that this notation of $\kappa$ is different from the notation given in [Asch1]).

We define the inverse mapping to $\kappa$ by $R:X \rightarrow K$. 
Let $\Phi$ be the homomorphism from $\la K \ra$ into $Sym(X)$. Then $\Phi(R(x)) = \rho_x$.

(f) For technical reasons it is useful to formally distinguish between elements of $K$ and elements of $X$,
as elements of $G$ may act on both sets, but in different ways.

Subloops,  homomorphisms, normal subloops, factor loops and simple loops are defined as usual in universal algebra:
A {\em subloop} is a  nonempty subset which is closed under loop multiplication. Be aware that we study finite loops.
Therefore, any subloops contains the identity.

{\em Homomorphisms} are maps between loops, which preserves loop multiplication. The preimage relation induces
an equivalence relation on the source loop, such that a product of equivalence classes is again an equivalence class.

{\em Normal subloops} are preimages of $1_{\circ}$ under a homomorphism and therefore subloops.
A normal subloop defines a partition of the loop into blocks (cosets), such that the set of products of elements from two blocks is again a block.
Such a construction gives factor loops as homomorphic images with the normal subloop as the kernel. 

{\em Simple loops} have only the full loop and the $1_{\circ}$-loop as normal subloops.

Finally we recall the definition of a {\em soluble loop} given in \cite{A}. A loop $X$ is {\em  soluble} if there
exists a series $1 = X_0 \leq \cdots \leq X_n= X$ of subloops with $X_i$ normal in $X_{i+1}$ and
$X_{i+1}/X_i$ an abelian group.\\

There are related concepts in the language of loop folders. We give here only the most important concept of a subfolder.
For other concepts and more results on loop folders see \cite{A} and \cite{AKP}.
\medskip\\
\noindent
{\bf Definition}
Let $(G,H,K)$ be a loop folder.
A {\em  subfolder} $(U,V,W)$ is a loop folder with $U \le G$, $V \le U \cap H$ and $W \subseteq U \cap K$.
\medskip\\

A subfolder defines a subloop $Y$ of a loop $X$, such that the multiplications $\circ_X$ and $\circ_Y$ coincide on $Y$.
Moreover, every subfolder is the folder of a subloop.

\begin{lemma}
\label{subfolders}
A subgroup $U \le G$ gives rise to a subfolder $(U,V,W)$, if and only if $U = (U \cap H)(U \cap K)$.
Then $V = U \cap H$ and $W = U \cap K$. 
\end{lemma}

\begin{bew}
Let $(U,V,W)$ be a subfolder. Then $W \leq U \cap K$ is a transversal to $V \le U \cap H$ in $U$,
which implies $U = WV \leq (U \cap H)(U \cap K)$. The Dedekind identity and the fact that $H \cap K = 1$
then implies $V = U \cap H$ and $W = U \cap K$. 

Now assume $U = (U \cap H)(U \cap K)$. Then $H \cap K = 1$ shows that $(U \cap K)$ is a transversal to
$(U \cap H)$ in $U$.
As $(U \cap K)$ acts transitively on the cosets of $(U \cap H)$ in $U$,
it follows that $(U \cap K)$ acts transitively on the cosets of $(U \cap H)^u$ in $U$ for every $u$ in $U$.
Thus $(U,U \cap H, U \cap K)$ is a subfolder.
\end{bew}

\begin{coro}
\label{HKsuper}
If $U$ is a subgroup of $G$ containing $H$ or $K$, then $(U, U \cap H, U \cap K)$ is a subfolder
of $(G,H,K)$.
\end{coro}

Though subfolders give access to inductive arguments, we have to be carefully  for two reasons.
\begin{itemize} 
\item A subfolder of a faithful loop folder may not be faithful.
\item A subfolder of a loop envelope may not be a loop envelope.
\end{itemize}

Another useful concept is the concept of morphisms between loop folders, see \cite{A}. 
We consider here only a special case, which is used to get faithful folders from arbitrary ones.

\begin{lemma}
\label{Hnormal}
Let $(G,H,K)$ be a loop folder to a loop $X$, $N \trianglelefteq G$ with $N \le H$ and $\overline{G}=G/N$.
Then $(\overline{G},\overline{H},\overline{K})$ is a loop folder to the same loop $X$.
\end{lemma}

\begin{bew}
The loop folder property is clearly inherited to the factor group.  
The two loops are natural isomorphic from the definition of the loop:  the multiplication depends only on the
action of $K$ on the $H$-cosets and $N$ is in the kernel of this action. 
\end{bew}

Finally we need another relation between $H=\Stab_G(1_{\circ})$ and $K=\{ \rho_x: x \in X \}$ in the Baer envelope of a loop.
Let $X$ be a loop and $x,y \in X$. Define: 
\[ h_{x,y} := \rho_x \rho_y \rho_{x \circ y}^{-1} \in \RMult(X). \]

\begin{lemma}
\label{Hgeneration}
Let $X$ be a loop and $(G,H,K)$ the Baer envelope of $X$. Then 
 \[ H= \langle h_{x,y}: x,y \in X \rangle. \]
\end{lemma}

\begin{bew} 
Let $H_1:=  \langle h_{x,y} : x ,y \in X \rangle \le G$. 
Then $H_1$ is a subgroup of $H$ and $|G:H|=|X|$. We claim that $G= H_1K$, which then yields the assertion.

We show this, using induction on the minimal length of elements $\sigma$ in $\RMult(X)$, expressed
as a product of elements of $K$. We assume that 
the minimal length is at least two, as words of length at most one are already in $K$.

Suppose $\sigma = \rho_{x_1} \rho_{x_2} \cdots \rho_{x_k} \in G$.  
If $k=2$, then $\sigma= h_{x_1,x_2} \rho_{x_1 \circ x_2} \in H_1 K$.  
For $k>2$, the word $\sigma_1:=\rho_{x_1 \circ x_2} \rho_{x_3} \cdots \rho_{x_k}$ has a shorter 
expression as product of elements from $K$, so $\sigma_1 = h_1 \rho_x$ for some $h_1 \in H_1$ and $x \in X$. 
Then $\sigma= h_{x_1,x_2} \sigma_1 = h_{x_1,x_2} h_1 \rho_x \in H_1 K$. 
\end{bew}

\subsection{Bol loops and twisted subgroups}

If we write the Bol identity using the right translations $\rho$, we get
\[  \mbox{ for all}~ y,z \in X: \quad \rho_{(y \circ z)\circ y} = \rho_y \rho_z \rho_y. \] 
This leads to the concept of twisted subgroups: 
\medskip\\
\noindent
{\bf Definition}
A {\em  twisted subgroup} $K$ of a group $G$ is a subset, such that for all $x,y \in K$: 
\begin{itemize}
\item[(1)] $1 \in K$
\item[(2)] $x^ {-1} \in K$ and 
\item[(3)] $xyx \in K$.
\end{itemize}

Notice that the second condition is not necessary for finite groups: If $x$ is in a twisted subgroup $K$,
then (1) and (3) imply that $K$ contains all the powers of $x$. Therefore, $K$ contains $x^{-1}$ as well.
 
We get:

\begin{lemma}
\label{Bol_twisted}
If $X$ is a Bol loop with faithful loop folder $(G,H,K)$, then $K$ is a twisted subgroup of $G$.
\end{lemma}

\begin{bew}
By Remark (d) we may assume that $(G,H,K)$ is the Baer envelope of the loop $X$. Therefore, the
elements of $K$ are the permutations $\rho_x$ with $x$ in $X$. Thus the Bol Identity implies that
$k_1k_2k_1$ is in  $K$ for all $k_1,k_2$ in $K$. 
\end{bew}

Notice, that for arbitrary (nonfaithful) loop folders $(G,H,K)$ to $X$, $K$ may not be a twisted subgroup:
we may replace elements $k$ of $K$ by $ck$ with $c \in$ Core$_G(H)$ without changing the loop multiplication.
\medskip\\
\noindent
{\bf Definition}
A loop folder $(G,H,K)$ to a Bol loop $X$ is called a {\em  Bol-folder}, if $K$ is a twisted subgroup of $G$.
\medskip\\

As just noted not any loop folder to a Bol loop is a Bol folder, but 
\begin{itemize}
\item The Baer folder of a Bol loop is a Bol folder, see \ref{Bol_twisted}.
\item Subfolders of Bol folders are Bol folders again.
\item Homomorphic images of Bol folders are Bol folders, see \ref{twisted_basic}(2).
\end{itemize}

In a Bol folder the Bol identity is translated into the group theoretic property of a twisted subgroup.
In this sense we call Bol folders  nice loop folders. 

We recall some of the results of Aschbacher on twisted subgroups from \cite{A1}.
As the original paper contains much more, we extract some of the critical arguments. 

\begin{lemma}
\label{twisted_basic}
Let $K$ be a twisted subgroup of the  group $G$. Then
\begin{enumerate}
\item[(1)] For all $k \in K$, $\langle k \rangle \le K$. 
\item[(2)] If $N \trianglelefteq G$, then the image of $K$ in $G/N$ is a twisted subgroup.
\item[(3)] For all $k \in K$, the set $k K$ is a twisted subgroup.  The twisted subgroups $k K$, $k \in K$ are called the associates of $K$.
\item[(4)] For $x \in K$, $xKx =K$.
\end{enumerate}
\end{lemma}

\begin{bew}
(1) is shown above. (2) is immediate from the definition. 
For (3) let $x,y,z \in K$. We write $(x y) (x z) (x y) = x (y x y)(y^{-1} z y^{-1}) (y x y)$ and $(x y)^{-1}=x (x^{-1} y^{-1} x^{-1})$. 
As $x^{-1} \in K$, it follows that  $1 \in xK$  and that $xK$ is a twisted subgroup of $G$. 
(4) follows from the definition.
\end{bew}

Let $G$ be a group with a twisted subgroup $K$ such that $G = \la K \ra$.
Define a sequence of relations $R_i \subseteq G \times G$ by
$$R_0 = \{(1,1)\}~\mbox{and}~R_{i+1} = \{(kx,k^{-1}y):(x,y) \in R_i, k \in K\}$$
and set $R_\infty = \bigcup\limits_{i=0}^\infty R_i$. As $G$ is finite, $R_\infty$ is a finite union of the $R_i`$s.
For $(x,y) \in R_\infty$ we write $x \equiv y$.

\begin{lemma}
\label{twisted}
Let $G, K$ and $\equiv$ as above, then the following holds.
\begin{enumerate}
\item[(1)] If $g_1 \equiv h_1$ and $g_2 \equiv h_2$, then $g_1 g_2 \equiv h_1 h_2$. 
\item[(2)] If $g \equiv h$, then $gK = Kh$.
\item[(3)] $\{ gK: g \in G \} = \{ Kg:g \in G \}.$
\item[(4)] $\Xi_K(G):=\{ g \in G: g \equiv 1 \}$ is a normal subgroup of $G$. 
\item[(5)] $\Xi_K(G) \subseteq \Psi_K(G):=\{ g \in G: gK =K \}$ is also a normal subgroup of $G$.
\item[(6)] $\Psi_K(G)K = K = K\Psi_K(G)$ and $\Psi_K(G) \subseteq K$
\item[(7)] If $\Xi_K(G)=1$, then there exists some $\tau \in \Aut(G)$ with $g\equiv g^\tau$ for all $g \in G$. 
Furthermore $\tau^2=1$, $k^\tau=k^{-1}$ for all $k \in K$ and the set $\Lambda:=\tau K\subset G\la \tau \ra$ is $G$-invariant. 
Notice, that by the action of $\tau$ on $K$ there is at most one automorphism of $G$ with that action. 
\end{enumerate}
\end{lemma}

\begin{bew}
(1) is obvious, as $g \equiv h$ if and only if there exist $k_1,...,k_n \in K$ with $g=k_1k_2 \cdots k_n$ and $h=k_1^{-1}k_2^{-1} \cdots k_n^{-1}$. 

In (2) use induction:  Let $(kg,k^{-1}h)$ be in $R_{i+1}$ with $k \in K$ and $(g,h) \in R_i$.
Then $gK = Kh$ and therefore $kgK = kKh = kKkk^{-1}h = Kh$ by \ref{twisted_basic}(4).

(3) is a consequence of (2): for any $g_1,h_1 \in G$ elements $g_2,h_2 \in K$ exists with $g_1 \equiv g_2$ and $h_2 \equiv h_1$. 

For (4) we use (1): If $g \equiv 1$, then $g^k=k^{-1} g k \equiv k 1 k^{-1} = 1$.    

In (5) notice, that $\Psi_K(G)$ is a subgroup containing $\Xi_K(G)$ by definition. Let $g \in \Psi_K(G)$ and $k \in K$.
Then $g^k K = k^{-1} g k K = k^{-1} (g (k K k)) k^{-1} = K$. As $G = \langle K \rangle$ we get (5).

For (6) notice, that $\Psi_K(G)K=K$ from the definition, so as $1 \in K$,  $\Psi_K(G) \subseteq K$. As  $kK= Kk^{-1}$ by (2), also $K \Psi_K(G) = K$. 

For (7) suppose $\Xi_K(G)=1$. Notice, that for any $g \in G$: if $h_1 \equiv g$ and $h_2 \equiv g$, then $h_1^{-1} h_2 \in \Xi_K(G)=1$, so $h_1=h_2$.
As $G=\langle K \rangle$, for any $g \in G$ there is a unique $h \in G$, such that $g \equiv h$. Define $g^\tau=h$ and notice, that $\tau$ is a
homomorphism by (1) with image $\langle K \rangle=G$, so an automorphism.  
As $K \subseteq C_G(\tau^2)$, but $\langle K \rangle =G$, $\tau^2=1$. Finally let $k \in K$.
Then $(\tau K)^k = k^{-1} \tau K k = \tau \tau k^{-1} \tau  K k = \tau (k K k) = \tau K$.     
\end{bew}

\noindent
{\bf Remarks.} (a) As $\tau$ acts on $\Psi_K(G)/\Xi_K(G)$ by inverting all elements, this section is abelian. 

(b) It may happen that $\tau$ is the identity. This happens for instance in Bol loops of exponent $2$.
\medskip\\

Following Aschbacher, $G$ is said to be {\em  reduced}, if $\Xi_K(G)=1$. 
Together with Lemma~\ref{twisted}(7) we get the following statement.

\begin{lemma}
\label{twisted_reduce}
Let $K$ be a twisted subgroup of the group $G$ and $G= \langle K\rangle$. 
Suppose, there exists an automorphism $\sigma \in \Aut(G)$ with $k^\sigma=k^{-1}$ for all $k \in K$. 
Then $\Xi_K(G)=1$ and $\sigma=\tau$. 
\end{lemma}

\begin{bew}
Let $g \in \Xi_K(G)$. There exist elements $k_1,...,k_n \in K$ with $g = k_1 k_2 \cdots k_n$ and $1 = k_1^{-1} k_2^{-1}  \cdots k_n^{-1}$. 
Using the automorphism property of $\sigma$ and its values on $K$, we get $\sigma(g)=1$, so $g=1$ and $\Xi_K(G)=1$.
Now we use  \ref{twisted}(7).
\end{bew}

This yields the following characterization of a twisted subgroup.

\begin{lemma}
\label{twisted_characterize}
Let $G$ be a finite group, $\tau \in \Aut(G)$ with $\tau^2=1$ and $K \subseteq G$ with $k^\tau=k^{-1}$ for all $k \in K$ and $\langle K \rangle=G$.
Then $K$ is a twisted subgroup, if and only if $1 \in K$ and $\Lambda = \tau K\subseteq G \langle \tau \rangle$ is $G$-invariant. 
\end{lemma}

\begin{bew}
If $K$ is a twisted subgroup, then we can use \ref{twisted_reduce} to get $\Xi_K(G)=1$. Now $\tau$ is the uniquely determined automorphism
defined in \ref{twisted}(7) and the statement holds. 

Suppose $\Lambda$ is $G$-invariant. Notice, that $(\tau k)^2= k^\tau k =1$. Let $k_1,k_2 \in K$.
Then $k_1  k_2  k_1 = \tau (\tau k_2)^{k_1} \in \tau \Lambda$. As $1 \in K$ and $k_1 \in K$ iterating this procedure, we see that $K$ contains all positive powers
of $k_1$. 
\end{bew}

\begin{lemma}
\label{normal_xi}[Asch1, 6.5]
Let $(G,H,K)$ be the Baer envelope to a Bol loop. 
Then $(\Xi_K(G), 1, \Xi_K(G))$ is a subfolder to a normal subloop  $\Xi(X)$.
Moreover, $\Xi(X)$ is a group and isomorphic to the group $\Xi_K(G)$.
\end{lemma}

Following Aschbacher, a Bol loop $X$ is called {\em  radical free}, if $\Xi(X)=1$. 

\subsection{$A_r$-loops}

If we wish to apply group theory in loop theory, the loops should have some automorphisms. 
Furthermore there should be a way to find other subloops than just those mentioned in \ref{subfolders}. 
A concept, which occurs naturally here is the concept of $A_r$-loops.

As it will turn out, Bruck loops are examples of $A_r$-loops, while general Bol loops need not have the $A_r$-property.
\medskip\\
\noindent
{\bf Definition}
The loop $X$ is called an $A_r$-loop, if for all $x,y \in X$: $ h_{x,y} \in \Aut(X)$, 
This means that 
\[ \mbox{ for all} ~x,y,u,v \in X: \quad (u \circ v)^{h_{x,y}} = (u)^h_{x,y} \circ (v)^h_{x,y} . \]

This definition implies, that subloops and homomorphic images of an $A_r$-loop are again $A_r$-loops.
Due to \ref{Hgeneration} our definition of $A_r$-loop is the same as in Section 4 of \cite{A}.
The following lemmata are results of Aschbacher, see Section 4 in \cite{A}.

\begin{lemma}
\label{Ar_criterion}[Asch1, 4.1]
A loop $X$ with Baer envelope $(G,H,K)$ is an $A_r$-loop if and only if $H$ acts on $K$ via conjugation. 
In that case $\rho_x^h = \rho_{x^h}$ for any $x \in X$, $h \in H \le Sym(X)$. 
\end{lemma}

\noindent
{\bf Definition}
An {\em  $A_r$-loop folder} is a loop folder $(G,H,K)$, such that $H$ acts on $K$ by conjugation. 
\medskip\\

In an $A_r$-folder the $A_r$-property of the loop is translated into the group theoretic condition, that $H$ acts by conjugation on $K$. Therefore, we call it a nice folder. It holds:

\begin{itemize}
\item The Baer folder of an $A_r$- loop is an $A_r$-folder, see \ref{Ar_criterion}. 
\item Subfolders of $A_r$-folders are $A_r$-folders again, see \cite{A} 4.2(2).
\item Homomorphic images of $A_r$-folders are $A_r$-folders, see \cite{A}  4.2(2).
\end{itemize}

We give an example of a loop folder to an $A_r$-loop, which is not an $A_r$-folder:
Let $G=\langle a,b|a^2=b^2=(ab)^4=1\rangle \cong D_8$,  $H:=\langle ab \rangle$ and $K:=\langle a \rangle$. 
The corresponding loop is the group of size 2.

The next lemma is essentially 4.3 of \cite{A}. 
\begin{lemma}
\label{ArFolders}[Asch1, 4.3]
Let $(G,H,K)$ be an $A_r$-loop folder to a loop $X$ and $L \le H$.  Then 
\begin{enumerate}
\item[(1)] $Fix_X(L):=\{ x \in X: x^l=x ~\mbox{ for all }~ l \in L \}$ is a subloop of $X$ and so closed under $\circ$. 
\item[(2)] For $k \in K$, $\{ h \in H: h^k \in H \} = C_H(k)$. 
\item[(3)] For $k \in K$, $\{ k \in K: L^k = L \} = \{ k \in K: [L,k]=1 \} =:C_K(L)$. 
\item[(4)] $(C_G(L),C_H(L),C_K(L))$ as well as $(N_G(L),N_H(L),C_K(L))$ are subfolders to $Fix_X(L)$. 
\item[(5)] $H$ controls $G$-fusion in $H$. 
\item[(6)] $[\langle K \rangle,\core_G(H)]=1$.
\end{enumerate}
\end{lemma}

Notice, that (6) is a consequence of (2).

If $(G,H,K)$ is an $A_r$-loop folder, then for $L \le H$, $N_G(L)$ and $C_G(L)$ give subfolders.
This is the reason, why the group theoretic approach to loops is so powerful:
The corresponding subloops may not be that interesting in loop theory, but the subgroups $C_G(L)$ and $N_G(L)$
play an important part in the local structure of a group.

\subsection{Bruck loops}
Recall, that a Bruck loop $X$ is a Bol loop, such that the inverse map $\iota: X \to X, x \mapsto x^{-1}$ is an automorphism of $(X,\circ)$.

\begin{lemma}
\label{bruck_baer_envelope}[Asch1, 6.6; AKP, 5.1]
Let $X$ be a Bruck loop with Baer envelope $(G,H,K)$. Then 
\begin{enumerate}
\item[(1)] $X$ is radical free.
\item[(2)] The map $\iota: X \to  X, x \mapsto x^{-1}$ induces on $G$ exactly $\tau$ from \ref{twisted}(7), 
\item[(3)] The set $\Lambda = \tau K \subseteq G\langle \tau \rangle$ is $G$-invariant.
\item[(4)] $H \le C_G(\tau)$. 
\item[(5)] $X$ is an $A_r$-loop.
\item[(6)] ${\rm Fix}_X(\iota)$ is a Bol loop of exponent 2 with folder $(C_G(\tau),H,C_K(\tau))$. 
\end{enumerate}
\end{lemma}

\begin{lemma}
\label{bruckloop}[AKP, 5.1]
Let $X$ be a Bol loop with Baer envelope $(G,H,K)$. The following statements are equivalent:
\begin{enumerate}
\item[(1)] $X$ is a Bruck loop.
\item[(2)] $X$ is a radical free $A_r$-loop.
\item[(3)] $\Xi(X)=\Xi(G)=1$ and $H$ acts on $K$ by conjugation.
\item[(4)] $\Xi(X)=\Xi(G)=1$ and $H \le C_G(\tau)$ for some $\tau \in \Aut(G)$ with $\tau^2=1$ and $k^\tau=k^{-1}$ for all $k \in K$.
\end{enumerate}
\end{lemma}

\noindent
{\bf Definition}
A {\em  Bruck folder} $(G,H,K)$ is a loop folder to a Bruck loop $X$, which is both an $A_r$-folder and a Bol folder.
So $K$ is a twisted subgroup and $H$ acts on $K$ by conjugation.
\medskip\\

Notice, that the Baer folder to a Bruck loop is a Bruck folder.

\begin{lemma}
\label{BruckFolder}
Let $(G,H,K)$ be a Bruck folder. Then the following hold.
\begin{enumerate}
\item[(1)]  There is a subgroup $Z$ of  $Z(\langle K \rangle)$ such that $\langle K \rangle /Z \cong \RMult(X)$.
\item[(2)]  There exists a unique $\tau \in \Aut(G)$ with $[H,\tau]=1$ and $k^\tau = k^{-1}$ for all $k \in K$. 
\item[(3)]  The set $\Lambda= \tau K \subseteq \Aut(G)$ is $G$-invariant.
\item[(4)]  Subfolders and homomorphic images are Bruck folder.
\end{enumerate}
\end{lemma}
\begin{bew}
(1) By \ref{ArFolders}(6), $\langle K \rangle$ is a central extension of $\RMult(X)$ by a group $Z \le H$
with $Z \le Z(\langle K \rangle)$. 

(2) If $\Xi_K(\la K \ra) = 1$, then $\tau$ exists by \ref{twisted}(7).
We claim that in fact $\Xi_K(\la K \ra) = 1$. By \ref{twisted}(4)-(6) $\Xi_K(\la K \ra) \subseteq K$.
Let $\alpha$ be the natural homomoprhism from $\la K \ra$ to $\RMult(X)$. Then $K$ and $\alpha(K)$ are
 twisted subgroups and $\alpha(\Xi_K(\la K \ra)) \leq \Xi_{\alpha(K)}(\alpha(\la K \ra)) = 1 $ by 
\ref{bruckloop}. Hence, $\Xi_K(\la K \ra) \leq$ ker $\alpha \leq Z(\la K \ra)$ by (1). As $|X| = |K| = |\alpha(K)|$,
we get $K \cap$ ker $\alpha = 1$. Thus $\Xi_K(\la K \ra) = 1$.

We claim that for all $h \in H \cap \la K\ra$ it holds $h^\tau = h$. We have $$k^{(\tau^h)}
= (k^h)^{\tau h^{-1}} = ((k^{-1})^h)^{h^{-1}} = k^{-1}.$$ Thus $\tau$ and $\tau^h$ are two automorphisms
which invert every element in $K$, which implies $\tau = \tau^h$ by \ref{twisted_reduce}.
So, $h^\tau = h$ for all $h \in H \cap \la K\ra$.

We extend the map $\tau$ to $G$ by defining $\tau(hk)= hk^{-1}$ for $h \in H$ and $k \in K$. Then $\tau$ is in $\Aut(G)$:

Let $h_1,h_2 \in H$  and $k_1,k_2 \in K$. Then
$$(h_1k_1h_1k_2)^\tau = (h_1h_2k_1^{h_2}k_2)^\tau = (h_1h_2k_1^{h_2}k_2k_3k_3^{-1})^\tau$$
where $k_3$ is an element in $K$ such that $k_1^{h_2}k_2k_3$ is an element in $H$. Then by the
definition of $\tau$ $$(h_1h_2k_1^{h_2}k_2k_3k_3^{-1})^\tau =h_1h_2k_1^{h_2}k_2k_3 k_3.$$
Notice also that $k_1^{h_2}k_2k_3 = (k_1^{h_2}k_2k_3)^\tau = (k_1^{h_2})^{-1}k_2^{-1}k_3^{-1}$.
This yields $$(h_1k_1h_1k_2)^\tau = h_1h_2(k_1^{h_2})^{-1}k_2^{-1}k_3^{-1}k_3 = h_1h_2(k_1^{h_2})^{-1}k_2^{-1}
= h_1 k_1^{-1} h_2k_2^{-1} $$$$= (h_1k_1)^\tau (h_1k_2)^\tau,$$ which yields the claim.
This gives (2).

As $K$ is a twisted subgroup (3) holds by \ref{twisted_characterize}.

Subfolders and homomorphic images of $A_r$-loop folders (resp. Bol folders) are again $A_r$-loop folders
(resp. Bol folders).
As subloops and homomorphic images of Bruck loops are again Bruck loops, we get (4).
\end{bew}

We add, that Bruck folders are nice in our sense:
The Bruck loop property of the loop (Bol identity and automorphic property)
translates into the existence of an automorphism $\tau \in \Aut(G)$ with
\begin{itemize}
\item $\tau^2=1$,
\item for all $h \in H, k \in K$: $h^\tau=h$ and $k^\tau=k^{-1}$.
\item $\tau K$ is $G$-invariant and 
\item $H$ acts by conjugation on $K$. 
\end{itemize}

Notice, that $\tau K$ is another transversal to $H$ in $\langle G,\tau \rangle$. If $\tau \ne 1$, $1 \not\in \tau K$, so in general
this transversal does not give a loop.
\medskip\\
\noindent
{\bf Notation}. Let $(G,H,K)$ be a Bruck folder $(G,H,K)$ and $\tau \in \Aut(G)$ the automorphism introduced in \ref{BruckFolder}(2). 
Then let \[{\bf  G^+} := G\langle \tau \rangle,\] the semidirect product of $G$ with $\tau$,
\[{\bf  H^+}:= H \langle \tau \rangle \le G^+~\mbox{and}~
{\bf  \Lambda }:= \tau K \subseteq G^+.\]  By \ref{BruckFolder}(3) $\Lambda$ is a $G^+$-invariant set of involutions.

\subsection{Bruck loops of 2-power exponent}
As mentioned in the introduction, results of Glauberman \cite{GG1}, \cite{GG2}, Aschbacher \cite{A}
and  Aschbacher, Kinyon and Phillips \cite{AKP}, now focus the attention to Bruck loops
of 2-power exponent. Again the loop properties translate into a property of $G$ and we get yet another nice folder type. 

\begin{lemma}
\label{BX2P_tau}
\cite{AKP},(5.13).
Let $(G,H,K)$ be a Bruck folder. 
Then $\tau \in O_2(G^+)$ if and only if every element of $K$ has $2$-power order. 
\end{lemma}

\begin{bew}
Suppose $\tau \in O_2(G^+)$. As $k^\tau = k^{-1}$ for all $k$ in $K$, it follows that every element in $K$
is of $2$-power order.

Now assume that every element in $K$ is of $2$-power order. Let $g = hk$ with $h $ in $H$ and $k$ in $K$.
Then $\tau^g = \tau ^k = k^{-2}\tau$. Hence, $\langle \tau,\tau^g \rangle$ is a 2-group for all $g \in G$. 
By the Baer-Suzuki Theorem $\tau$ is in $O_2(G^+)$. 
\end{bew}

Notice, that if $(G,H,K)$ is a Bruck folder to a loop $X$, 
$X$ is of exponent 2 iff $K$ is a union of $1\in G$ and $G$-conjugacy classes of involutions.
\medskip\\
\noindent
{\bf Definition}
Let $X$ be a Bruck loop of 2-power exponent.
A loop folder $(G,H,K)$ to $X$ is called a {\em  BX2P-folder},
if it is a Bruck folder and every element of $K$ has a $2$-power order.
Equivalently $\tau \in O_2(G^+)$. 
A loop folder $(G,H,K)$ is called a {\em  BX2P}-envelope, if $(G,H,K)$ is a BX2P-folder and a loop envelope, so $G = \langle K \rangle$. 
\medskip\\

Again, the Baer folder to a Bruck loop of 2-power exponent is a BX2P-folder, while
subfolder and images of BX2P-folders are again BX2P-folders. 

\begin{lemma}
\label{overlineK}
Let $(G,H,K)$ be a BX2P-folder. Then $k^2 \in O_2(G)$ for all $k$ in $K$. 
If $\overline{G}=G/O_2(G)$, then $1 \in \overline{K}$ and $\overline{K}-\{1\}$ is a union of conjugacy classes of involutions of $\overline{G}$.
\end{lemma}

\begin{bew}
Let $k \in K$ and $\tau \in \Aut(G)$ be the automorphism of \ref{BruckFolder}(2). 
Then $k^2 = \tau k^{-1} \tau k = [\tau,k] \in [O_2(G^+),G] \le O_2(G^+) \cap G \le O_2(G) $ by \ref{BX2P_tau}.
In particular $\overline{K}=\overline{\Lambda}$ in $G^+/O_2(G^+)$. As $\Lambda$ is a union of $G^+$-conjugacy classes
of involutions (and $1$ if $\tau = 1$) by \ref{BruckFolder}(3), also the last statement holds. 
\end{bew}

\begin{rem}
The set $K$ has not to be a normal set in $G$, but $\Lambda$ is a normal set in $G^+$ by \ref{BruckFolder}(3).
If $\tau = 1$, then $K$ is normal in $G$. As $\la K \ra = G \cap \la \Lambda\ra$, the group $\la K \ra$ is normal in
$G$.
\end{rem}

This lemma is the reason, why the special case of Bol loops of exponent 2 is so closely related to the general case
of Bruck loops of 2-power exponent. 

While working on the case of Bol loops of exponent 2, we decided to completely ignore the structure of $O_2(G)$,
as almost nothing was known about $O_2(G)$.
Luckily in the general case of Bruck loops of 2-power exponent, the group $G^+/O_2(G^+)$ behaves exactly
as in the special case of Bol loops of exponent 2, since the sets $K$ and $\Lambda$ have the same image modulo $O_2(G^+)$:
$\overline{K}= \overline{\Lambda}$. 
This trick was already used in \cite{AKP} to reuse the arguments of \cite{A} for the classification of $M$-loops.
\section{The Bruck loops of 2-power exponent}
This section contains just about anything, which was previously known about Bruck loops of 2-power exponent,
as well as about Bol loops of exponent 2. We formulate this knowledge in the language of BX2P-folders.
Not everything is needed in order to prove the main theorem, but most statements are used in [S] and in [BS3].

\subsection{Basic results}
Here we present results known before the paper \cite{A} of Aschbacher. 
Most of the arguments essentially go back to Heiss, see \cite{Hei}.

In a Bol loop, the order of every element divides the loop order. Therefore, the following holds.

\begin{lemma}
\label{evensize}
\begin{itemize}
\item[(1)] A Bruck loop of 2-power exponent has even size or size 1. 
\item[(1)] If $(G,H,K)$ is a BX2P-folder, then $|G:H|=|X|$ is 1 or even. 
\end{itemize}
\end{lemma}
\begin{bew}
For $k \in K$, $\langle k \rangle$ acts semiregularly on the $H$-cosets of $G$ by the loop folder
property. 
\end{bew}

\begin{lemma}
\label{Hsuper}
Let $(G,H,K)$ be a BX2P-folder. Let $U \le G$ be a subgroup such that 
 $U = (U \cap H)(U \cap K)$. Then  the subfolder  to $U$  (see \ref{subfolders}) is $(U,U\cap H,U \cap K)$ and the size of 
the corresponding subloop is $|U:U\cap H|$. 
In particular overgroups of $H$ and of $\langle K \rangle$ satisfy this condition.
\end{lemma}

\begin{bew} See \ref{subfolders}. \end{bew}

\begin{lemma}
\label{noHinvert}
Let $(G,H,K)$ be a BX2P-folder. Let $\lambda \in \Lambda$, $h \in H$ and $g \in G$.
If $(h^g)^\lambda = (h^g)^{-1}$, then $h^2 = 1$. 
\end{lemma}

\begin{bew}
Suppose $h^{g\lambda}=(h^g)^\lambda = (h^g)^{-1} = (h^{-1})^g$. Let $\mu=g \lambda g^{-1} \in \Lambda$. Then $h^\mu=h^{-1}$,
so $[h,\mu]=h^{-2} \in H$. But $[h,\mu] = h^{-1} \mu h \mu = \mu^h \mu \in \Lambda \Lambda \cap H = K K \cap H$. Since $KK \cap H=1$ by the loop folder
property, $h^2=1$.   
\end{bew}

\begin{lemma}
\label{O2prime}
Let $(G,H,K)$ be a BX2P-folder. Then
the following holds.
\begin{itemize}
 \item[(1)] $O_{2'}(G) \le C_H(\langle K \rangle)$.
\item[(2)] If $(G,H,K)$ is a faithful BX2P-folder, then $O_{2'}(G)=1$.
\item[(3)] If $(G,H,K)$ is a BX2P-envelope, then $O_{2'}(G) \le Z(G) \cap H$.
\end{itemize}
\end{lemma}

\begin{bew}
$O_{2'}(G) H$ gives rise to a subfolder by \ref{HKsuper}, but $|O_{2'}(G) H : H|$ is odd, so by \ref{evensize},
$|O_{2'}(G)H:H|=1$. By \ref{noHinvert} then $[\langle \Lambda \rangle, O_{2'}(G)]=1$. 
\end{bew}

The following stronger version of \ref{evensize} holds. It has very strong consequences, see 
Corollaries~\ref{2Himplies2G}, \ref{O22prime}.

\begin{lemma}
\label{even_index}
Let $(G,H,K)$ be a BX2P-folder and $U \le G$ with $H \le U$. Then $|G:U|$ is even or $1$.
\end{lemma}

\begin{bew}
Assume $|G:U|$ to be odd. Then $U$ contains a Sylow-2-subgroup of $G$, so every element of $K$ is conjugate
to some element of $U \cap K$. Let $U^+=U \langle \tau \rangle \le G^+$, so $|G^+:U^+|$ is odd.
Then $|\{ \lambda^g: \lambda \in \Lambda\cap U^+, g \in G \}| \le 1+(|U^+:H^+|-1)|G^+:U^+| = 1 + |G:H| - |G:U|$.
Since $|G:H|=|K|=|\Lambda| = |\{ \lambda^g: \lambda \in \Lambda\cap U^+, g \in G \}|$ we get $|G:U|=1$. \end{bew}

\begin{coro}
\label{2Himplies2G}
$H$ is a $2$-group if and only if $G$ is a $2$-group.
\end{coro}

\begin{bew}
If $H$ is a 2-group, then $H$ is a contained in 2-Sylow $M$ of $G$, so by \ref{even_index} $|G:M|=1$
and $G$ is a 2-group.
\end{bew}

\begin{coro}
\label{O22prime}
Let $(G,H,K)$ be a BX2P-folder. Then $O_{2,2'}(G)H = O_2(G) H$.  
\end{coro}

\begin{bew}
$O_2(G) H$ is of odd index in $O_{2,2'}(G) H$, so the statement is a consequence of \ref{even_index}
and \ref{Hsuper}.
\end{bew}

\begin{lemma}
\label{soluble_loops}
Let $(G,H,K)$ be a  BX2P-envelope to a soluble Bruck loop $L$ of $2$-power exponent.
Then $|L|=|G:H|$ is a power of $2$.
\end{lemma}

\begin{bew}
As $L$ is soluble, there is a series $1 = L_0 \leq \cdots \leq L_n = L$ of subloops of $L$ with $L_i$
normal in $L_{i+1}$ such that $L_{i+1}/L_i$ is an abelian group. Every element in $L$ is of $2$-power,
which implies that $L_{i+1}/L_i$ is a $2$-group. Thus $|L|$ is a power of $2$.
\end{bew}

In this case even more can be said.

\begin{lemma}
\label{2powerorder}
Let $(G,H,K)$ be a BX2P-envelope and $|G:H|$ a power of $2$. Then $G$ is a $2$-group.
\end{lemma}

\begin{bew}
As $|G:H|$ is a 2-power, $H$ contains Sylow subgroups for all odd primes.
But then the product of any two elements of $K$ has to be of 2-power order: 
If $k_1 k_2 $ is not of $2$-power order for $k_1,k_2 \in K$, then $\tau k_1^{-1} \tau k_2 \in \Lambda \Lambda$
is not of 2-power order. Then there exist $\lambda_1,\lambda_2 \in \Lambda$ with $\lambda_1 \lambda_2$ of odd prime
order and $\lambda_1$ inverts $\lambda_1 \lambda_2$. By \ref{noHinvert} this is a contradiction, as $\lambda_1 \lambda_2$ is conjugate
into $H$ by assumption. Now by the Baer-Suzuki Theorem, $\langle \Lambda \rangle $ is a 2-group, so $G = \langle K \rangle$ is a 2-group too.  
\end{bew}

Now we study conditions on the BX2P-envelope which force $L$ to be soluble.

\begin{coro}
\label{solubleloop}
Let $(G,H,K)$ be a BX2P-envelope to a Bruck loop $L$ and $G = O_2(G)H$. Then $L$ is soluble.	
\end{coro}
\begin{bew}
By \ref{2powerorder} $G$ is a $2$-group. Let $G = G_r \unrhd G_{r-1} \unrhd \cdots \unrhd G_0 =H$ be a normal series
starting at $G = G_0$ and ending at $H$ such that $G_{i+1}/G_{i}$ is of order $2$. Let
$L_i$ be the loop defined by $(G_i, H, G_i \cap K)$. Then $|L_{i+1}- L_i| = |L_i|$. This property allows to
construct an homomorphism from $L_{i+1}$ into $\ZZ_2$ with kernel $L_i$.
Thus $L_{i+1}$ is a normal subloop of $L_i$, for
$0\leq i \leq r-1$, and $L_{i+1}/L_i$ is a group of order $2$. This shows that $L = L_r$ a soluble.
\end{bew}

Lemmas~\ref{soluble_loops} and \ref{2powerorder} imply: if a Bruck loop of $2$-power exponent is soluble,
then $G$ is a $2$-group. The following theorem shows that if $G$ is soluble, then the Bruck loop is
soluble as well, and $\la K \ra$ is already a $2$-group.

\begin{thm}
\label{solublegroups}
Let $(G,H,K)$ be a BX2P-folder and assume that $G$ is soluble. Then $\langle K \rangle \le O_2(G)$ is a $2$-group. 
\end{thm}

\begin{bew}
Our proof uses \ref{inv_invert}, which was introduced in \cite{A}. Nevertheless the Theorem was already proved in \cite{Hei}.

Let $\overline{G}=G/O_2(G)$. By \ref{O22prime}, $F^\ast(\overline{G})= F(\overline{G}) \le \overline{H}$. 
Let $\lambda \in \Lambda$. If $\lambda$ acts nontrivially on $F(\overline{G})$, it inverts some element of odd prime order $p$ 
in $F(\overline{G})$. By \ref{inv_invert}, $\lambda$ inverts some element of order $p$ in the preimage of $F(\overline{G})$,
but $H$ contains a Sylow-$p$-subgroup of that preimage. By \ref{noHinvert} we get a contradiction. 
Therefore, the elements in $\Lambda$ act trivially on $F(\overline{G})$. As $C_{\overline{G^+}}(F(\overline{G^+})) \le Z(F(\overline{G^+}))$,
it follows $\lambda \in O_2(G)$, which implies $\langle \Lambda \rangle \le O_2(G^+)$. 
\end{bew}

\begin{coro}
\label{Characsoluble}
Let $(G,H,K)$ be a BX2P-folder to a Bruck loop $L$. Then
$L$ is soluble if and only if $\la K\ra$ is a $2$-group.
\end{coro}

The following lemma will be helpful.

\begin{lemma}
\label{oddnormal}
Let $(G,H,K)$ be a BX2P-folder and $\overline{G}=G/O_{2'}(G)$.
Then $(\overline{G},\overline{H},\overline{K})$ is a loop folder to the same loop.
\end{lemma}

\begin{bew}
By \ref{O2prime}, $O_{2'}(G) \le H$, so \ref{Hnormal} gives the result.
\end{bew}

\subsection{Selected results of Aschbacher, Kinyon, Phillips}

Next we present some of the results from \cite{A}, and \cite{AKP}, which are fundamental to our results.

For the next lemma see also \ref{ArFolders} and \ref{Hsuper}.

\begin{lemma} 
\label{subloops}
Let $(G,H,K)$ be a BX2P-folder to a Bruck loop $X$ of $2$-power exponent.
\begin{itemize}
\item[(1)]  Let $L \le H$. Then $(N_G(L),N_H(L), C_K(L))$ and $(C_G(L),C_H(L),C_K(L))$ are BX2P-subfolders to a (the same) subloop of $X$. 
\item[(2)]  Let $U \le G$ with $|U| \le | U \cap H| |U \cap K|$. Then $(U,U \cap H,U \cap K)$ is a BX2P-subfolder to a subloop of $X$. 
\end{itemize}
\end{lemma}

\begin{bew}
By \ref{BruckFolder}(4), subfolder of BX2P-folders are BX2P-folders. So
(1) is a consequence of \ref{ArFolders}(4) and 
(2) follows from \ref{Hsuper}.
\end{bew} 

The idea to ignore $O_2(G)$ resp. $O_2(G^+)$ origins in \cite{A}. We present here the main arguments:

\begin{lemma}
\label{inv_invert}
(\cite{A},8.1(1))
Let $G$ be a group and $x \in G$ an involution.
If $\overline{x} \in \overline{G}:=G/O_2(G)$ inverts some element $\overline{y} \in \overline{G}$ of odd prime order $p$, then $x$ inverts some element of order $p$ in $G$.
\end{lemma}

Now we get the next lemma, which will be used repeatedly in [S].

\begin{lemma}
\label{noHoverlineinvert}
Let $(G,H,K)$ be a BX2P-folder. Let $\overline{G}= G/O_2(G)$ and $\overline{x} \in \overline{K}, \overline{y}\in \overline{G}$.
If $1\neq o(\overline{y})$ is odd and $\overline{y}^{\overline{x}} = \overline{y}^{-1}$, 
then for every $\overline{z} \in \overline{G}$: $\overline{y}^{\overline{z}} \notin \overline{H}$. 
In particular, $\overline{y} \notin \overline{H}$.   
\end{lemma}

\begin{bew}
Assume otherwise. Let $\tau \in \Aut(G)$ be the automorphism defined above and recall that  $\tau k \in \Lambda$ 
and that $O_2(G^+) k = O_2(G^+) \tau k$, as $\tau \in O_2(G^+)$. 
Since $\langle \overline{y}, \overline{x} \rangle$ is a dihedral group with all involutions
conjugate, we may assume w.l.o.g that $o(\overline{y})$ is some odd prime $p$, by replacing $\overline{y}$ with some suitable element from
$\langle \overline{y} \rangle$.

We can choose preimages $x \in K$ of $\overline{x}$ and $y \in H$ of $\overline{y}$ with $o(y)=o(\overline{y})$. Recall,
that 
$\tau x$ is another preimage of $\overline{x}$ in $G^+$. As $\tau x$ inverts some element of prime order $p$ in 
$$\overline{ O_2(G^+)\langle y\rangle },$$  by \ref{inv_invert} then $\tau x$ inverts some element of prime order $p$ in $O_2(G) \langle y \rangle$. 
But $O_2(G) \langle y \rangle \le O_2(G) H$ and $H$ contains a $p$-Sylow-subgroup of $O_2(G)H$. 
So $\tau x$ inverts some element of odd order, which is conjugate into $H$, a contradiction to \ref{noHinvert}. 
\end{bew}

The following definition is taken from \cite{AKP}.
\medskip\\
\noindent
{\bf Definition}
\label{Mloop}
An $M$-loop is a finite Bruck loop $X$, such that each proper section of $X$ is soluble, but $X$ itself is not soluble. 
\medskip\\

The next theorem is Theorem~3 of \cite{AKP}.

\begin{thm}
\label{AMT} \cite{AKP}
Let $X$ be an $M$-loop with Baer envelope $(G,H,K)$, $J = O_2(G)$ and $G^\ast = G/J$. Then
\begin{enumerate}
\item[(1)] $X$ is a Bruck loop of 2-power exponent.
\item[(2)] $G^\ast \cong \PGL_2(q)$ with $q = 2^n+1 \ge 5$, $H^\ast$ is a Borel subgroup of $G^\ast$ and $K^\ast$ consists of the involutions in $G^\ast - F^\ast (G^\ast)$.
\item[(3)] $F^\ast(G) = J$.
\item[(4)] Let $n_0 = |K \cap J|$ and $n_1 = |K \cap a J|$ for $a \in K -J$. Then $n_0$ is a power of $2$, $n_0 = n_1 2^{n-1}$ and $|X|=|K| = (q+1) n_0 = n_1 2^n (2^{n-1}+1)$. 
\end{enumerate} 
\end{thm}

The following lemma is another formulation of Aschbachers \cite{A} (12.5)(2), which is based on an idea of Heiss.
The formula for Bruck loops occurs in (3.2)(3) of \cite{AKP} and will be heavily used in [S] and [BS3].

\begin{lemma}
\label{HeissEquation}
Let $(G,H,K)$ be a BX2P-folder and $N \trianglelefteq G$ with $O_2(G) \le N$. 
Let $a_i, i \in \{1,...,r\}$ be representatives for the orbits of $\overline{G}= G/N$ on $\overline{\Lambda}^\sharp$,
$m_i:=| \{ \overline{a_i}^{\overline{G}} \} |$, $n_i= |\Lambda \cap a_i N^+|$ and $n_0:=\Lambda \cap N^+$. Then 
\[ |K| = | \Lambda| = n_0 + \sum_{i=1}^r n_i m_i. \] 
\end{lemma}

\begin{bew}
Let $\Lambda_i := \{ a \in \Lambda : \overline{a} \in \overline{a_i}^{\overline{G}} \}$ and $\Lambda_0 := \Lambda \cap N^+$.
Then $\{ \Lambda_i: i \in \{0,..,r\} \}$ is a partition of $\Lambda$ with $|\Lambda_0|=n_0$ and $|\Lambda_i|=n_i m_i$ for $i \in \{1,..,r\}$. 
\end{bew}
 
\subsection{Additional results}

In the following $\overline{G}$ always denotes the group $G/O_2(G)$.
The results presented here emerged during work on the classification of passive simple groups.
We start with a corollary to \ref{HeissEquation} which is basic to the classification of the passive
groups in [S]. It is a very powerful tool to get a full Sylow $p$-subgroup of $G$ into $H$ for $p$ a prime divisor of $|H|$.

\begin{coro}
\label{HeissPrime}
Let $(G,H,K)$ be a  BX2P-folder. 
Suppose $O_2(\overline{H})=1$  and that there exists an odd prime $p$ dividing $|\overline{G}|$, such that
$m_i \equiv 0 \pmod{p}$ for all $i \in \{1,..,r\}$, with $m_i$ as in \ref{HeissEquation} for $N=O_2(G)$. 
Then $p$ does not divide $|K|=|G:H|$.
\end{coro}

\begin{bew}
Since by assumption $O_2(H) \subseteq O_2(G)$, we have $O_2(O_2(G)H) = O_2(G)$. Therefore, as $|O_2(G)H:H|$ is a 2-power,
$$(O_2(G)H,H,O_2(G)H \cap K)$$ is a subfolder to a soluble subloop by \ref{solubleloop}.
Hence $|O_2(G)H \cap K|$ as well as $|\langle O_2(G)H \cap K\rangle |$ is a 2-power, and, as $\langle O_2(G)H \cap K\rangle$
is normal in $O_2(G)H $, we get
$ \langle O_2(G) H \cap K \rangle \le O_2(O_2(G)H) = O_2(G)$, 
Thus
$n_0=|O_2(G)^+ \cap \Lambda|=|O_2(G) \cap K|=|O_2(G)H \cap K|$ is a 2-power. By \ref{HeissEquation}  $p$ does not divide $|K|$.  
\end{bew}

There is a corollary to \ref{noHoverlineinvert}, which generalizes Theorem \ref{solublegroups}:

\begin{coro}
\label{ZeroComponentCase}
Let $(G,H,K)$ be a BX2P-folder. If $F^\ast(\overline{G})= F(\overline{G})$, then $\overline{G}=\overline{H}$.
\end{coro}

\begin{bew}
We have $F(\overline{G}) \le \overline{H}$ by \ref{O22prime}. By \ref{noHoverlineinvert}, no element of $\overline{K}$
acts nontrivially on $F(\overline{G})$. Therefore $\langle \overline{K} \rangle \le C_{\overline{G}}(F(\overline{G})) \le Z(F(\overline{G}))
\le \overline{H}$.
Therefore,  $\langle \overline{K} \rangle = 1$ and $\overline{G}= \overline{H}$.  
\end{bew}

Therefore, by \ref{solubleloop}. and \ref{ZeroComponentCase}, in a nonsolvable loop with BX2P-folder $(G,H,K)$, 
$\overline{G}$ has components. 
The following lemma makes use of soluble subloops.
It shows, that $\overline{H}$ has to contain certain elements of odd order.

\begin{lemma}
\label{O_upper_2_criterion}
Let $(G,H,K)$ be a BX2P-folder, $\overline{G}=G/O_2(G)$ and $U\le G$ be a subgroup with the following properties:\\
\begin{itemize}
\item[(1)] $U = (U \cap H)(U \cap K)$.
\item[(2)] $[O_2(U), O^2(U)] \le O_2(G)$.  
\item[(3)] $\langle U \cap K \rangle \le O_2(U)$.
\end{itemize}
Then $O^2(U) \le O_2(G)H$. 
\end{lemma}

\begin{bew}
Let $u \in U$ be of odd order. We can write $u = h k$ with $h \in H \cap U$ and $k \in K \cap U$ by (1).
Now $k \in \langle K \cap U \rangle \le O_2(U)$ by (3). By (2)  we have $[u,k] \in [O^2(U),O_2(U)] \le O_2(G)$. 
In $\overline{G}=G/O_2(G)$ the element $\overline{u}$ is of odd order and commutes with $\overline{k}$. 
As $[\overline{u}, \overline{k}] =1$ implies $[\overline{h}, \overline{k}] =1$ and as $\overline{k}$ is of order
$1$ or $2$, it follows that $\overline{k} \in \overline{H}$, which yields the assertion.
\end{bew}

There exists a powerful generalization to nonsoluble subloops.

\begin{lemma}
\label{passive_centralizing_components}
Let $(G,H,K)$ be a BX2P-folder, $\overline{G}=G/O_2(G)$ and $D:=\langle K \rangle$.
Then $O^2(C_{\overline{G}}(\overline{D})) \le \overline{H}$. 
\end{lemma}

\begin{bew}
Let $x \in G$ be of odd order, such that $[\overline{D},\overline{x}]=1$. We can write $x= h k$ with $h \in H, k \in K$. 
As $k \in O_2(G) D$, $[\overline{k},\overline{x}]=1$. As $\overline{x}=\overline{h} \overline{k}$, $[\overline{h},\overline{k}]=1$,
so $\overline{k}$ is in $\overline{H}$ as $\overline{x}$ has odd order. Therefore $\overline{x}$ is in $\overline{H}$. 
\end{bew}

\noindent
{\bf Definition}
A Bruck loop $L$ of 2-power exponent is called a {\em  2M-loop}, if $L$ is not soluble, but every proper
subloop is soluble.

\begin{rem}
Notice, that an $M$-loop has to be simple while a $2M$-loop may not. For instance, a nonsplit extension of a soluble subloop
by a simple non-soluble loop may be a $2M$-loop. In order not to have to exclude such extensions, we have
introduced the concept of a $2M$-loop. 
\end{rem}

The classification of $M$-loops by Aschbacher, Kinyon and Phillips given in Theorem~\ref{AMT} yields a 
description of the $2M$-loops.

\begin{lemma}
\label{qLemma}
Let $q>1$ be an integer with $q-1$ a 2-power. Then $q=2$ or $9$ or $q \geq 5$ is a Fermat prime.
\end{lemma}

\begin{bew} See \cite{BScommuting} for a proof, based on Zsigmondy's Theorem. \end{bew}

\begin{lemma}
\label{2NloopEnvelope}
Let $(G,H,K)$ be a BX2P-envelope to a $2M$-loop $L$. 
Then the following holds.
\begin{itemize}
\item[(1)] $C_G(O_2(G)) \le O_2(G)$,
\item[(2)]$\overline{G} \cong \PGL_2(q)$ and $q=9$ or $q\ge 5$ is a Fermat prime,
\item[(3)] $|G:O_2(G)H|=q+1$,
\item[(4)] $\overline{K}$ consists of 1 and all involutions in $\PGL_2(q)\setminus{\PSL_2(q)}$,
\item[(5)] $O_2(G)=(O_2(G)\cap H)(O_2(G) \cap K)$.
\end{itemize}  
\end{lemma}

\begin{bew}
Let $L_1$, $L_2$ be normal proper subloops. These subloops are soluble by definition of the $2M$-loop.
Notice, that $L_1 L_2$ is another soluble normal subloop, see [Bruck], thus a proper subloop too.
Therefore there exists a biggest proper normal subloop $L_0$, which is soluble. The quotient $L/L_0$
then is an $M$-loop as defined in \ref{Mloop}. Let $D:=\langle R(x): x \in L_0 \rangle \le G$. 
Then, as $L_0$ is a normal soluble subloop, $D$ is a normal $2$-group of $G$, so $D \le O_2(G)$ and $G/D$ is a loop envelope to an $M$-loop, see [Asch] 2.6. If we manage to prove the statement 
for $(\tilde{G},\tilde{H},\tilde{K})$ with $\tilde{G}= G/D$, the statement holds for $(G,H,K)$,
so we may assume $D=1$. 

The structure of a faithful loop envelope to an $M$-loop is described in Theorem \ref{AMT}, which together with \ref{qLemma} implies the statement. Notice, that (5) follows from Theorem~\ref{AMT}(4).

Now assume that $(G,H,K)$ is not faithful. By \ref{ArFolders}(6) $C:= \core_G(H) $ is in $Z(G)$. 
Let $Z:=O_{2'}(Z(G))$. Then $Z\le C$ by \ref{O2prime}(1) and
$(\tilde{G},\tilde{H},\tilde{K})$ with $\tilde{G}:= G/C$ is a faithful loop envelope to an $M$-loop
by \ref{Hnormal}. So we can apply Theorem \ref{AMT}. Then $\overline{G}= G/O_2(G)$
is a central extension of $\PGL_2(q)$ with $\overline{Z}$ still contained in the group generated by 
$\overline{K}$. Thus, if $Z \cong \overline{Z} \neq 1$, then 
 $q=9$ and $|Z|=3$, as this is the only case of nontrivial odd order Schur multiplier of the groups in question.
(The $r$-part of the Schur multiplier of a perfect group may be nontrivial for noncyclic Sylow-$r$-subgroups only. 
The unique noncyclic case $q=9$ actually results in a Schur multiplier $\ZZ_3$ for $\Alt(6) = \PSL_2(9)$.)

According to Theorem \ref{AMT}
the involutions in $\overline{G}\setminus{\overline{G}^\prime}$ are in $\overline{K}$. However  in this case, involutions outside $\overline{G}^\prime$ invert $Z$,
as is visible using the embedding of $3\Alt(6)$ into $\SL_3(4)$, see \cite{ATLAS} p.23 for the action
of $L_3(4)$-automorphisms on the Schur multiplier.
This contradicts \ref{O22prime} and \ref{noHoverlineinvert}, so $Z=1$.

The factorization $O_2(G)=(O_2(G)\cap H)(O_2(G) \cap K)$ can be seen as follows:
We have $O_2(G) H = H (O_2(G)H \cap K)$ by \ref{Hsuper}. Let $k \in K \cap O_2(G)H$.
As $\overline{H}\cong q:(q-1)$ does not contain involutions of $\PGL_2(q) \setminus{\PSL_2(q)}$, we obtain that $\overline{k}=1$.
Thus $k \in O_2(G)$ and the assertion follows with the Dedekind identity. 
\end{bew}

A powerful application of Aschbachers results is the $2M$-loop-embedding: Any nonsoluble Bruck loop of 2-power exponent
contains a 2M-loop as a subloop. Since the structure of a $2M$-loop is very restricted, we get strong information on $G$. 

\begin{lemma}
\label{2NloopEmbedding}
Let $(G,H,K)$ be a BX2P-folder with $G \ne O_2(G) H$. 
Then some subgroup $U \le G$ exists such that
\begin{itemize}
\item $U = (U\cap K)(U\cap H), U = \langle U \cap K \rangle$,
\item The loop to $(U,U\cap H,U \cap K)$ is a $2M$-loop,
\item $F^\ast(U) = O_2(U)$, 
\item $U/O_2(U) \cong \PGL_2(q)$ for $q\ge 5$ a Fermat prime or $q=9$,
\item $|U:O_2(U) (U\cap H)|=q+1$,
\item $\overline{K\cap U}$ consists of 1 and all involutions in $\PGL_2(q)\setminus{
\PSL_2(q)}$,
\item There exist elements of order $\frac{q+1}{2}$ in $U$ inverted by elements of $\Lambda\cap U^+$.
\item There exist elements $h \in U \cap H \cap G^{(\infty)}$ of order $3$ or $q$ if  $q=9$ or $q\neq 9$, 
respectively.
\item In particular $G^{(\infty)}$ contains a section isomorphic to $\PSL_2(q)$.  
\end{itemize}
\end{lemma}

\begin{bew}
We can find the subgroup $U$ recursively: If the loop is nonsoluble,
but every subloop is soluble, the loop is itself a  $2M$-loop. 
Else we can find a proper nonsoluble subloop, which contains a $2M$-loop $Y$. 
Set $U:= \la R(x): x \in Y\ra$. Then $(U, U\cap H, \cap K)$ is a loop folder to $Y$, see [Asch1] 2.1.

Now \ref{2NloopEnvelope} describes the structure of $U$, which implies the statements.
\end{bew}

\section{The Proof of Theorem~\ref{EnvelopeGroups}}

If not explicitely defined otherwise, $\overline G= G/O_2(G)$ and for subsets $X \subseteq G$, $\overline{X}$ is the image
of the natural homomorphism from $G$ onto $\overline{G}$.
 \medskip\\
\noindent
{\bf Definition}
\label{anchor_def}
Let $S$ be a finite non-abelian simple group. Let ${\cal L}_S$ be the class of all Bruck loops $X$ of 2-power exponent,
such that there is  a BX2P-folder $(G_X,H_X,K_X)$ to $X$ with $F^\ast(G_X/O_2(G_X)) \cong S$.
A prime $p, p>2$, is called {\em  passive} for $S$, if $p \nmid |X|$ for all $X \in {\cal L}_S$. 
($p$ itself may not divide $|S|$.)

The smallest passive prime $p \in \pi(S)$ is called the {\em  anchor prime of S}.
It is the smallest odd prime $p \in \pi(S)$ such that for every $X \in {\cal L}_S$,
$p$ does not divide $|G_X:H_X|=|X|$.

The finite non-abelian simple group $S$ is called {\em  passive}, if 
every odd prime $p \in \pi(S)$ is passive. 
\medskip\\

\noindent
{\bf Remarks}
(1) Notice that the definition of a passive finite non-abelian simple group is equivalent to the definition
given in the introduction:

The finite non-abelian simple group $S$ is passive if and only if $X$ is soluble for every $X \in {\cal L}_S$
if and only if $G = O_2(G) H$ whenever $(G,H,K)$ is a BX2P-folder with $F^\ast(\overline{G}) \cong S$.
The equivalence of these conditions follows from \ref{solubleloop}, the $2M$-loop embedding
\ref{2NloopEmbedding} and \ref{soluble_loops}, \ref{2powerorder}: The $2M$-loop embedding implies, that $\overline{G}^{(\infty)}$ contains 
elements of order either $3$ or $5$, which are products of two elements in $\overline{K}=\overline{\Lambda}$,
so any $2M$-loop embedding prevents one of the primes $3$ or $5$ from being passive. 

(2) The anchor prime to a finite non-abelian simple group may not exist. Its existence will be established 
later by classifiying the non-passive finite simple groups, using the classification of finite simple groups.

(3) If $S$ is passive, then $S$ has an anchor prime, usually 3, except in case of the Suzuki groups ${}^2B_2(q)$, where it is 5.

\begin{lemma} 
Let $S\cong \PSL_2(q)$ for $q\ge 5$ a Fermat prime. Then either $q$ or $3$ is the anchor prime of $S$. 
\end{lemma}

\begin{bew}
The 2M-loop embedding, \ref{2NloopEmbedding}, and the list of subgroups of $\PSL_2(q)$ by Dickson
and the fact that $5$ does not divide $q+1= 2^{4n}+1+1$, implies that we have an embedding such that $UO_2(G) = G$ 
and $U \cong \PGL_2(q)$.

We get that $H$ always contains a Sylow-$q$-subgroup of $U$. Thus the prime $q$ is passive for $S$. For $q=5$ the existence of examples ensures, that $q=5$ is the smallest such prime.
In the other cases there may be no examples of $M$-loops for the corresponding $q$, so $\PSL_2(q)$
is passive. Then $q=3$ is the anchor prime. If examples exist, the anchor prime is $q$.  
\end{bew}

\begin{lemma}\label{AnchorA6}
Let $S \cong \PSL_2(9) \cong \Alt(6)$. Then $p=3$ is the anchor prime.
\end{lemma}

\begin{bew}
Let $(G,H,K)$ be a BX2P-folder with $F^\ast(\overline{G}) \cong S$. If $G=O_2(G)H$, then $H$ contains a Sylow-3-subgroup of $G$.
By \ref{2NloopEmbedding} and Dixons theorem we can only embedd
2M-loops for $q=5$ or $q=9$. The case $q=9$ implies, that $H$ contains  a Sylow-3-subgroup of $G$.

Otherwise there is a subgroup $U$ in $G$ such that $U/O_2(U) \cong PGL_2(5)$ and such that 
$\overline{U} \cap \overline{H} \cong 5:4$ by \ref{2NloopEmbedding}. Then
$H$ contains elements of order 5. These elements are inverted
by inner involutions of $\Alt(6)$ and (if $\overline{G}$ contains  $\PGL_2(9)$) involutions of $\PGL_2(9)$
outside $\PSL_2(9)$. Therefore, by \ref{noHoverlineinvert}  $\overline{K}$ can consist only of the 1-element, the 15 transpositions of $\Sym(6)$
and the 15 involutions of $\Sym(6)$, which are a product of three commuting transpositions.
Therefore, $|\overline{G}:\overline{H}| \leq 31$. As $\overline{G}$ is a subgroup of $Aut(\Alt(6))$, it follows
from the list of subgroups of $Aut(\Alt(6))$ that $\overline{H} \cong Sym(5)$.
Thus $H$ contains an element $x$  of order 3. Then $(C_G(x),C_H(x),C_K(x))$ is a subfolder by \ref{subloops}(1).
As $C_G(x)$ contains a Sylow-3-subgroup of $G$ which covers $O_{2^\prime}\overline{(C_G(x))}$, 
the subgroup $H$ contains a Sylow-3-subgroup of $G$ by Lemma~\ref{O22prime}.
Thus $3$ is the anchor prime to $\Alt(6)$.
\end{bew}

\noindent
{\bf Definition}
Let $(G,H,K)$ be a BX2P-folder and $C$ a component of $\overline{G}=G/O_2(G)$. 
An {\em  anchor group} $A$ of $C$ is a subgroup of $C \cap \overline{H}$ such that $A \in \Syl_p(C)$ for the anchor prime $p$ of $C/Z(C)$.
\medskip\\

The following proposition is crucial for the proof of Theorem~\ref{EnvelopeGroups} as it will be used to show that every component of $\overline{G}$ is either normal in 
$\langle \overline{K} \rangle$ or contained in $\overline{H}$. 

The assumption, that every simple section has an anchor prime can be considered as a kind of ${\cal K}$-group assumption:
In the classification of finite simple groups, ${\cal K}$ is the list of `known' finite simple groups and the goal was to show, 
that ${\cal K}$ contains every finite simple group. 

With regard to Bruck loops we first study groups, such that every simple section has an anchor prime.
In [S] it is shown that every finite simple group has an anchor prime.

\begin{prop} 
\label{AnchorGroupsExist}
Let $(G,H,K)$ be a BX2P-folder and suppose every non-abelian simple section of $G$ has an anchor prime. 
Then every component of $\overline{G}$ has an anchor group. 
\end{prop}

\begin{bew}
The proof proceeds by induction on $|G|$. We reduce the structure of $\overline{G}$ in multiple steps and produce either anchor groups
or a contradiction. 
\medskip\\
(1) $O_{2'}(G)=1$:\\
If $O_{2'}(G)\neq 1$, then by induction on $G/O_{2'}(G)$, the statement holds for the loop folder from \ref{oddnormal}. 
Since $O_{2'}(G) \le H$ by \ref{O2prime}, the statement holds in $G$ too.
\medskip\\
(2) $F(\overline{G})=1$: \\
By \ref{O22prime} we have $F(\overline{G}) \le \overline{H}$. If $\overline{x}\in F(\overline{G})$
for some element $x \in H$ of odd prime order, then $(C_G(x),C_H(x),C_K(x))$ is a subfolder by \ref{subloops}(1).
Since $O_{2'}(G) = 1$ by (1), $C_G(x)$ is a proper subgroup.
Let $E$ be the full preimage of $C_{\overline{G}}(\overline{x})$ in $G$. Then by Frattini
$E = O_2(G)C_G(x)$ which yields that 
$C_G(x)$ covers $C_{\overline{G}}(\overline{x})$. Clearly,  the latter contains $E(\overline{G})$.
Therefore anchor groups of components of $C_G(x)/O_2(C_G(x))$, which exist by induction, lift to anchor groups
of $\overline{G}$. 
\medskip\\
(3) $E(\overline{G})$ contains more than one component: \\
Else $\overline{G}$ has a unique component, which has an anchor prime $p$ by assumption. By definition
of the anchor prime  an anchor group exists. \\
\medskip\\
(4) If $C \cap \overline{H}$ contains non-trivial elements of odd order for some component $C$ of $\overline{G}$,
then anchor groups for all components exists: \\
Let $x$ be such an element. Then $C_G(x)$ covers all but the component $C$. By induction we get anchor groups
for all components of $C_G(x)/O_2(C_G(x))$. These lift to anchor groups for the components of $\overline{G}$,
other than $C$. Since we have more than one component, we can use some element $z$ of odd prime order
in one of these anchor groups to get the anchor group of $C$ by induction on $C_G(z)$, which covers $C$.
\medskip\\
(5) $\overline{H} \cap E(\overline{G})$ is a 2-group: \\
Otherwise let $\overline{x} \in \overline{H} \cap E(\overline{G})$ be of odd prime order $p$. 
We can write $\overline{x}$ uniquely as $\overline{x} = \overline{x_1} \overline{x_2} \cdots \overline{x_k}$ with $\overline{x_i} \in C_i$, where $C_1,...,C_k$ are the components
of $\overline{G}$. 

If $\overline{x_i} = 1$ for some $i$, $C_G(x)$ covers the component $C_i$, so by induction on $C_G(x)$ we get an anchor group to $C_i$ as in (4).
We saw already in (4), that this implies, that all components have anchor groups.
So $\overline{x_i} \ne 1$ for every $i$.
Now $C_{E(\overline{G})}(\overline{x})$ is the direct product of the $C_{C_i}(\overline{x_i})$. 
In particular $\langle \overline{x_1},\overline{x_2},...,\overline{x_k} \rangle \le O_p(C_{E(\overline{G})}(\overline{x})) \le O_p(C_{\overline{G}}(\overline{x}))$.
Let $x$ be some preimage of $\overline{x}$ of order $p$. 

Since $C_G(x)$ covers $C_{\overline{G}}(\overline{x})$, it follows that $O_p(C_{\overline{G}}(\overline{x}))$ is covered
by $O_{2,2'}(C_G(x))$. By \ref{O22prime}, we may choose therefore preimages of the $\overline{x_i}$ in $H$. 
By (4) we now get anchor primes for all components of $\overline{G}$.

By  \ref{2NloopEmbedding} there is an element $h \in H$  of odd prime order $p$. 
\medskip\\
(6)  $\overline{h}$ normalizes every component of $\overline{G}$: \\
Otherwise let $C$ be a component with $C^{\overline{h}} \ne C$ and $D = CC^{\overline{h}} \cdots C^{\overline{h}^{p-1}}$, 
the closure of $C$ under $\overline{h}$. Now $C_D(\overline{h}) = \{ c c^{\overline{h}} \cdots c^{\overline{h}^{p-1}}: c \in C \} \cong C$.

By \ref{subloops}(1),  $C_G(h)$ is a group to a subloop. Notice, that $C_D(\overline{h})$ maps to a component of $C_G(h)/O_2(C_G(h))$:
$D$ is subnormal in $\overline{G}$, so $C_D(\overline{h})$ is subnormal in $C_{\overline{G}}(\overline{h})$,
but $C_G(h)$ covers $C_{\overline{G}}(\overline{h})$. 

By induction, we get an anchor group $A$ of $C_D(\overline{h})$. But then $A \le E(\overline{G}) \cap \overline{H}$, so
$E(\overline{G}) \cap \overline{H}$ contains elements of odd order contrary to (5).
 \medskip\\
(7) We get anchor groups for all components of $\overline{G}$:\\
We use \ref{2NloopEmbedding} to get an additional property of $h \in H$: 
some $h \in H$ of odd prime order exists, such that $h \in N_h^{(\infty)}$, with $N_h$ the normal closure of $h$. (Recall, that the element $h$ is in a $\PSL_2(q)$-section.)

Let $G_1$ be the subgroup of $G$ consisting of all elements, which normalize every component of $\overline{G}$.
Notice, that the preimage $E$ of $E(\overline{G})$ is contained in $G_1$. By the Schreier-conjecture
$G_1/E$ is soluble. By (6) we have $h \in G_1$. Therefore $N_h \le G_1$. As $h \in N_h^{(\infty)} \le E$, this is a contradiction to (5). 
\end{bew} 

The following lemma reveals the idea behind the anchor groups: Anchor groups insure that  the involutions of $\overline{K}$  fix all components of $\overline{G}$ (see \ref{overlineK}).

\begin{lemma}
\label{AnchorPreventsWreath}
Let $(G,H,K)$ be a BX2P-folder and suppose that every non-abelian simple section of $G$
has an anchor prime.
Then every element $x$ of $K$ normalizes every component $C$ of $\overline{G}$.
In particular a component of $\overline{G}$ is either a normal subgroup of  $\langle \overline{K} \rangle$ or contained in $\overline{H}$. 
\end{lemma}

\begin{bew}
Let $x \in K$, $\lambda \in \Lambda$ with $\overline{x}=\overline{\lambda}$ and $C$ be a component of $\overline{G}$. 
Assume $C^x \ne C$. Let $A,B$ be anchor groups to the components
$C$ and $C^x$, respectively, which exist by \ref{AnchorGroupsExist}. 
As $C$ and $C^x$ are isomorphic, the corresponding anchor primes $p_1$ and $p_2$ are equal.

In particular $A B \in \Syl_{p_1}(C C^x)$. Let $y \in A$ be of order $p_1$. 
As $p_1$ is odd and $A$ is Sylow in $C$, not every element of order $p_1$ of $A$ is in $Z(C)\ge C \cap C^x$.
Therefore, we may choose $y \notin C^x$.

Then $\overline{x}$ inverts  the element $y^{-1} y^{\overline{x}}$, which is of order $p_1$, and hence is conjugate to some element of $AB \le \overline{H}$.
This is a contradiction to \ref{noHoverlineinvert}. 

So $[C, \langle \overline{K} \rangle] \le C \cap \langle \overline{K} \rangle$.
Therefore either $C \trianglelefteq \langle \overline{K} \rangle$ or $[C,\langle \overline{K} \rangle]=1$.
In the latter case let $c \in C$ be of odd order. We can write $c=\overline{k} \overline{h}$ 
with $k \in K$, $h \in H$. As $\overline{k}$ commutes with $c$, it follows that $\overline{k}$ commutes with
$\overline{h} = c\overline{k}$ as well. The fact that $c$ is of odd order and $\overline{k}$ an involution
yields $\overline{k}$ is contained in $\overline{H}$, which implies that $c$ is in $ \overline{H}$.
Now $C= O^2(C)$ yields $C \le \overline{H}$.
\end{bew}

{\bf Proof of Theorem~\ref{EnvelopeGroups}}
Let $(G,H,K)$ be a BX2P-envelope and assume, that every non-abelian simple section
of $G$ is either passive or isomorphic to $\PSL_2(q)$ for $q=9$ or a Fermat prime $q \ge 5$.  
If $\oG = \overline{H}$, then by \ref{2powerorder} $G$ is a $2$-group and the theorem holds.
Hence we may assume $\oG \neq \overline{H}$.
If $F^\ast(\oG) = F(\oG)$, then by \ref{ZeroComponentCase} $\oG = \overline{H}$, so we assume $F^\ast(\oG) \neq F(\oG)$.

We prove the theorem by induction on the order of $G$. 
\medskip\\
(1)$F(\overline{G}) \leq Z(\langle \overline{K} \rangle)$:\\
As $\overline{G}=\langle \overline{K} \rangle$ and as no element of $\overline{K}$ acts nontrivially on $F(\overline{G})$ by \ref{noHoverlineinvert} and \ref{O22prime}, $F(\overline{G}) \le Z(\langle \overline{K} \rangle)$.
\medskip\\
Recall that by \ref{AnchorPreventsWreath} 

\noindent
(2) every component of $\overline{G}$ is  normal in $\overline{G}$.
\medskip\\
(3) Every passive component $C$ which is not isomorphic to $\PSL_2(q)$, with $q = 9$ or $q \geq 5$ a Fermat prime 
is contained in $\overline{H}$:\\
 We distinguish the two cases that $\overline{G}$
contains either one or more components.

$\overline{G}$ has precisely one component $C$.  If $F(\overline{G}) = 1$, then by the definition of passive $C$ is as desired. Thus we may assume $F(\overline{G}) \neq  1$. Let $F$ be a subgroup of $H$ of odd order
 such that $\overline{F} = F(\overline{G})$.  Let $N$ be the full preimage of $C$ in $G$. Then, as 
$C\leq C_{\overline{G}}(\oF)$, $N=O_2(G) C_N(F)$ by Frattini. Let $G_1 := C_G(F)$, so $(G_1,G_1 \cap H,G_1 \cap K)$
is a subfolder by \ref{subloops}. Set $G_2 = G_1/F$. Then $(G_2 , \beta(H \cap G_1), \beta(K \cap G_1))$ with $\beta$ the natural homomorphism from $G_1$ onto $G_1/F$ is a loop folder to the same loop by \ref{oddnormal}.
As $\oF \leq Z(\oG)$ it follows $F(G_2/O_2(G_2)) = 1$.
Then $F^\ast(G_2/O_2(G_2) \cong C/Z(C)$. Hence by induction and by
the definition of a passive group $G_2 = O_2(G_2)\beta(H \cap G_1)$.
Therefore, the loop to the folder for $G_2$, as well as this one for $G_1$, is soluble by \ref{solubleloop}. Now
 \ref{soluble_loops} implies
$G_1 = O_2(G_1)(H\cap G_1)$. Hence $N$ is contained in $O_2(G)H$ and $C$ in 
$\overline{H}$.

$\oG$ has more then one component. Let $D$ be a component of $\oG$ different from $C$. By \ref{AnchorGroupsExist}
there is a nontrivial element $x$ in $H$ of odd order such that $\overline{x}$ is in an anchor group of $D$.
Then, as $\overline{x} $ is in $D$,  $[C,\overline{x}] = 1$. Let $ N$ be the full preimage of $C$ in $G$. Then by the Dedekind identity
$N = O_2(G)C_N(x)$. Let $G_1 = C_G(x)$, $H_1= C_H(x)$ and $K_1 = G_1 \cap K$. Then $(G_1,H_1, K_1)$ is a proper
subfolder of $(G,H,K)$ by \ref{subloops}. Then \ref{AnchorPreventsWreath} implies that $C\cong \overline{C_N(x)}$ is contained
in either  $\overline{\la K_1\ra} $ or in $\overline{H_1}$. In the first case we obtain by induction on $|G|$ ($|G_1| < |G|$) the  statement of Theorem~\ref{EnvelopeGroups} for $(\la K_1\ra,\la K_1\ra \cap H_1,K_1)$. Hence, $C \cong \PSL_2(q)$, with $q = 9$ or $q\geq 5$ a Fermat prime in contradiction
to our assumption. In the latter case $C \leq \overline{H_1}$ which yields $C \leq \overline{H}$, the assertion.
\medskip\\
(4) $\overline{H}$ does not contain a component of $\overline{G}$:\\
Assume $\overline{H}$ contains a component $C$ of $\overline{G}$. Let $\overline{x}$ be an element in $K$.
By (2) $[\overline{x},C] \leq C$. Set $B_x := C:\la x \ra$. By  \ref{noHoverlineinvert} $\overline{x}$
does not invert an element of odd order in $C$. Hence $\la \overline{x},\overline{x}^{\overline{b}}\ra$
is a $2$-group for every  $\overline{b}$ of $B$. Thus by Baer-Suzuki $\overline{x} \in O_2(B_x)$
and therefore $[C,\overline{x} ] = 1$. This implies $C \leq Z(\overline{G})$, which is not possible.
\medskip\\

(3) and (4) imply that\\
(5) Every component of $\oG$ is isomorphic to $\PSL_2(q)$, with $q = 9$ or $q\geq 5$ a Fermat prime. 
\medskip\\
(6) $\oG/F^\ast(\oG)$ is an elementary abelian $2$-group: \\
As $F(\oG) \leq Z(\oG)$ by (1), we have 
$\oG/Z(F^\ast(\oG))$ is isomorphic to a subgroup of $Aut(E(\oG))$ which fixes every component of $\oG$.
By (5) the outer automorphism group of every component of $\oG$ is an elementary abelian $2$-group which yields the
assertion.
\medskip\\ 
(7) If $\oG$ has a unique component, then the assertion holds:\\
By \ref{2NloopEmbedding} there is a subgroup $U$ of $G$ such that $E(\overline{U}) \cong \PSL_2(q^\prime)$,
$q^\prime = 9$ or $q \geq 5$ a Fermat prime and $U/O_2(U) \cong \PGL_2(q^\prime)$. By assumption, (5)
and by (33.14) of [Asch0] $E(\oG) \cong \PSL_2(q)$ or $3\PSL_2(9)$.

We claim that $q = q^\prime$. If $q \neq q^\prime$, then $q = 9$ and $q^\prime= 5$ by the subgroup list of 
$\PSL_2(q)$ given by Dickson and by the fact that $q+1 = 2^{2^n}+2 \equiv 3 (5)$ for $n >1$. 
Then by \ref{AnchorA6} $\overline{H}$ contains a Sylow $3$-subgroup of $\oG$. On the other hand
$\overline{H \cap U}$ also contains a Sylow $5$-subgroup by \ref{2NloopEmbedding}. Hence 
$\overline{H} = \overline{U}$ or $\oG$ which is not possible. Thus $q = q^\prime$.

This yields that $O_2(G)U = G$ or $q = 9$ and $\oG$ contains a subgroup isomorphic to $\PGL_2(9)$. 
So, it remains to consider the case $q = 9$. As $H$ contains a Sylow $3$-subgroup of $G$, \ref{noHoverlineinvert}
implies that $\overline{K}$ consists only of involutions in $\PGL_2(9)\setminus{\PSL _2(9)}$. 
As $\la K \ra = G$ we get $\oG \cong \PGL_2(9)$. If $\oG \cong 3\PGL_2(9)$, then the elements in 
$\overline{K}$ invert $Z(\oG)$ by [Atlas] p. 23,  in contradiction to  \ref{noHoverlineinvert} and \ref{O22prime}
 and (1) and (2) of the theorem. By \ref{2NloopEmbedding} (3) and (4) hold as well.
\medskip\\ 
(8)  If $\oG$ has at least two  components $C_1$ and $C_2$, then the assertion holds:\\
By \ref{AnchorGroupsExist} we get anchor groups $A_i \le C_i$.
Let $B_i \le H$ be of odd order with $\overline{B_i} = A_i$. We can use induction on $G_i:= \langle C_G(B_i) \cap K \rangle$ by applying \ref{subloops}(1). This shows that $\overline{G_i}$ is as described in the statement of the theorem. In particular, no components of $\overline{G_i}$ are isomorphic to $3\PSL_2(9)$ and the 
elements in $\overline{K} \cap \overline{G_i}$ induce $\PGL_2(9)$-involutions on these components.
As before, we see that $\overline{K}$ does not contain an element which induces an $\Sym(6)$-involution
on some component, so (1) and (2) hold. 

Moreover by induction every component of $\oG_i$ acts faithfully on $O_2(G_i)$.  As $O_2(G_i)$ is contained
in $O_2(G)$, it follows that $O_2(G) = F^\ast(G)$, which is (4). 
By induction and as $q(q-1)$ is a maximal subgroup of $\PSL_2(q)$
we get (3) using (4).
This proves the assertion.
\qed

\end{document}